\documentclass{amsart}

\usepackage{amssymb}
\usepackage{graphicx}
\usepackage{lineno}

\newtheorem{theorem}{Theorem}

\newcommand{\U}{{\mathcal U}}

\newcommand{\IC}{{\mathbb C}}

\newcommand{\ID}{{\mathbb D}}

\def\be{\begin{equation}}
\def\ee{\end{equation}}
\newcommand{\bthm}{\begin{theorem}}
\newcommand{\ethm}{\end{theorem}}
\newcommand{\beqq}{\begin{eqnarray*}}
\newcommand{\eeqq}{\end{eqnarray*}}

\begin{document}

\title[On certain properties of some subclasses of univalent functions]{On certain properties of some subclasses of univalent functions}

\author[M. Obradovi\'{c}]{Milutin Obradovi\'{c}}
\address{Department of Mathematics,
Faculty of Civil Engineering, University of Belgrade,
Bulevar Kralja Aleksandra 73, 11000, Belgrade, Serbia}
\email{obrad@grf.bg.ac.rs}

\author[N. Tuneski]{Nikola Tuneski}
\address{Department of Mathematics and Informatics, Faculty of Mechanical Engineering, Ss. Cyril and Methodius
University in Skopje, Karpo\v{s} II b.b., 1000 Skopje, Republic of North Macedonia.}
\email{nikola.tuneski@mf.edu.mk}



\begin{abstract}
In this paper we determine the disks $|z|<r\le1$ where for different classes of univalent functions, we have the property
$${\rm Re}\left\{2\frac{zf'(z)}{f(z)}-\frac{z f''(z)}{f'(z)}\right\}>0\qquad (|z|<r).$$
\end{abstract}

\maketitle

Let ${\mathcal A}$ denote the family of all analytic functions
in the unit disk $\ID := \{ z\in \IC:\, |z| < 1 \}$   satisfying the normalization
$f(0)=0= f'(0)-1$.

\medskip

Further, let $\mathcal{S}$ be the subclass of $\mathcal{A}$ consisting of all univalent functions in $\ID$, and  $\mathcal{S}^{\star}$ and $\mathcal{K}$ be the subclasses of
${\mathcal A}$ of functions that  are starlike and convex in $\ID$, respectively.
Next,  let $\mathcal{U} $ denote the set of all  $f\in {\mathcal A}$ satisfying the condition
$$\left |\left (\frac{z}{f(z)} \right )^{2}f'(z)-1\right | < 1 \qquad (z\in \ID).$$
More on this class can be found in \cite{OP-01,OP_2011,DTV-book}.

\medskip

Next, by $\mathcal{G}$ we denote the class of all  $f\in {\mathcal A}$ in $\ID$ satisfying the condition
$${\rm Re}\left\{1+ \frac{z f''(z)}{f'(z)}\right\}<\frac{3}{2}\qquad (z\in \ID).$$
More about the class $\mathcal{G}$ one can find in \cite{JO_95} and \cite{OPW_2013}.

\medskip
 In their paper (\cite{MM-73}) Miller and Mocanu introduced the classes of $\alpha$-convex functions
 $f\in \mathcal{A}$ by the next condition:
 \be\label{eq1}
{\rm Re}\left\{(1-\alpha)\frac{zf'(z)}{f(z)}+ \alpha \left(1+ \frac{z f''(z)}{f'(z)}\right)\right\}>0\qquad (z\in \ID),
\ee
where $\frac{f(z)f'(z)}{z}\neq 0$ for all $z \in \ID $, and $\alpha\in\mathbb{R}$.
Those classes they denoted by $\mathcal{M_{\alpha}}$ and proved the next

\medskip

\noindent
\bf {Theorem A}. \\
{\it
\indent $(a)$ $\mathcal{M_{\alpha}}\subseteq \mathcal{S}^{\star}$ for every $\alpha \in\mathbb{R}$;\\
\indent $(b)$  $\mathcal{M}_{1}=\mathcal{K}\subseteq \mathcal{M_{\alpha}}\subseteq \mathcal{S}^{\star}$ for $0\leq \alpha \leq 1$; \\
\indent $(c)$  $\mathcal{M_{\alpha}}\subset\mathcal{M}_{1}=\mathcal{K}$ for $\alpha > 1$}.

\medskip
\rm
In \cite{OT_2019}  the authors proved

\noindent
\bf {Theorem B}. {\it
\\
\indent $(a)$  $\mathcal{M_{\alpha}}\subset\mathcal{U}$ for $\alpha \leq -1$;\\
\indent $(b)$  $\mathcal{M_{\alpha}}$ is not subset of $\mathcal{U}$ for $0\le\alpha \leq 1$.
}

\medskip
\rm
Choosing $\alpha=-1$ in Theorem A(a) and Theorem B(a),  from \eqref{eq1}, we have that the condition
\begin{equation}\label{prop-1}
{\rm Re}\left\{2\frac{zf'(z)}{f(z)}-\frac{z f''(z)}{f'(z)}\right\}>1\qquad (z\in \ID)
\end{equation}
implies $f\in \mathcal{S}^{\star}\cap \mathcal{U}$, i.e., the above inequality is sufficient for
univalence in the unit disc. As expected, it is not necessary condition for univalence, i.e., univalent functions does not necessarily have property \eqref{prop-1}. See functions $f_2$ and $f_3$ analysed bellow.

\medskip
But, is the following weaker inequality necessary for univalence
$${\rm Re}\left\{2\frac{zf'(z)}{f(z)}-\frac{z f''(z)}{f'(z)}\right\}>0\qquad (z\in \ID)?$$
The answer is also negative. Even more, it is not necessary condition even for starlikeness, nor for the classes $\U$ and $\mathcal{G}$.

\medskip
Namely, let consider the  differential operator
\be\label{eq2}
D(f;z):=2\frac{zf'(z)}{f(z)}-\frac{z f''(z)}{f'(z)}
\ee
and the functions $k(z)=\frac{z}{(1-z)^{2}}$, $f_{1}(z)=\frac{z}{1-z^{2}}$, $f_{2}(z)=-\log(1-z)$, and $f_3(z)=\frac{z(1-\frac{1}{\sqrt{2}}z)}{1-z^{2}}$. Then, we have, respectively,
\[
\begin{split}
D(k;z)&=1+\frac{1+z^{2}}{1-z^{2}},\\
D(f_{1};z)&=1+\frac{1-z^{2}}{1+z^{2}},\\
D(f_{2};z)&=-\frac{z(2+\log(1-z))}{(1-z)\log(1-z)},\\
D(f_{3};z)&=\frac{-\sqrt{2}z^{3}+3z^{2}-3\sqrt{2}z+2}{(1-\frac{1}{\sqrt{2}}z)(1-\sqrt{2}z+z^{2})}.
\end{split}
\]

\medskip

From the previous remark we easily conclude that the functions $k$ and $f_{1}$ belong to the class
$\mathcal{S}^{\star}\cap \mathcal{U}$, but for the function $f_{2}$ (which is  convex) for
$z=r$, $0\leq r<1$, we have
$$D(f_{2};r)=-\frac{r(2+\log(1-r))}{(1-r)\log(1-r)}<0$$
if $ 1-e^{-2}=0.86466\ldots\leq r<1.$
Also, we note that  $f_{2}\notin\mathcal{U}$.

\medskip

For the function $f_3$, in \cite{PoOb-2005}, the authors showed that it is close-to-convex and univalent in $\ID$, but not in $\mathcal{U}$. Additionally, ${\rm Re}[D(f;z)]>0$  does not hold on the unit disk. Indeed, let we put
\be\label{eq9}
D(f_3;z)=:\frac{g(z)}{h(z)},
\ee
where
$$g(z)=-\sqrt{2}z^{3}+3z^{2}-3\sqrt{2}z+2$$
and
$$h(z)=\left(1-\frac{1}{\sqrt{2}}z\right)(1-\sqrt{2}z+z^{2}),$$
and use  $z=r$, $0\leq r<1$. Then it is evident that $ h(r)>0$ for all  $0\leq r<1$. Also we have
$g'(r)=-3(\sqrt{2}r^{2}-2r+\sqrt{2})<0 $ for all $r\in [0,1)$, which implies that $g(r)$ is a decreasing function
on the interval $[0,1)$. Thus, $2=g(0)\geq g(r)> g(1/\sqrt{2})=0$ for $0\leq r< 1/\sqrt{2}$, and
$g(r)\leq 0$ for $\frac{1}{\sqrt{2}}\leq r<1.$
Now, from \eqref{eq9}, we easily conclude that the condition ${\rm Re}[D(f;z)]>0$ is not satisfied for the function $f_3$
in the disc $|z|<r$, where $\frac{1}{\sqrt{2}}\leq r<1$, i.e., for close-to-convex functions, ${\rm Re}[D(f;z)]>0$ on a disk with radius smaller then $\frac{1}{\sqrt{2}}=0.7071\ldots.$

\medskip
The above analysis raises the question of finding radius $r_*$ for each of the classes defined above, such that ${\rm Re}[D(f;z)]>0$
at least in the disc $|z|<r_*$. The next theorem answers this question. We don't know if the values for $r_*$ are the best possible.

\medskip

\bthm\label{th 1}
Let $D(f;z)$ be  defined by  \eqref{eq2}. Then
$${\rm Re}[D(f;z)]>0\qquad (|z|<r_*)$$
in each of the following cases:
\begin{itemize}
  \item [(\emph{i})]  $f\in \mathcal{U}$ and  $r_*=r_1=0.839\ldots$ is the root of the equation
$r^{3}+2r^{2}-2=0$;
  \item [(\emph{ii})] $f\in \mathcal{S}^{\star}(1/2)$ and $r_*=r_2=\sqrt{\frac{\sqrt{5}-1}{2}}=078615\ldots$;
  \item [(\emph{iii})] $f\in \mathcal{G}$ and $r_*=r_3=\frac{2}{3}=0.666\ldots$;
  \item [(\emph{iv})] $f\in \mathcal{S}^{\star}$ and $r_*=r_4=\frac{1}{2}=0.5$;
  \item [(\emph{v})] $f\in \mathcal{S}$ and $r_*=r_5=\frac{1}{4}=0.25.$
\end{itemize}
\ethm

\begin{proof}

  (i)  First, from the definition of the class $\mathcal{U}$, we easily conclude that
$f\in\mathcal{U}$ if, and only if, there exists a function $\phi$, analytic in $\ID$ with
$|\phi(z)|\leq 1$ in $\ID$, such that
\be\label{eq3}
\left[\frac{z}{f(z)} \right]^{2}f'(z)=1+z^{2}\phi(z).
\ee
From \eqref{eq3}, after some calculations, we obtain that
\be\label{eq4}
2\frac{z f'(z)}{f(z)}-\frac{z f''(z)}{f'(z)}=2\frac{1-\frac{1}{2}z^{3}\phi'(z)}{1+z^{2}\phi(z)}.
\ee
Since $|\phi(z)|\leq 1$, then
\be\label{eq5}
|\phi'(z)|\leq \frac{1-|\phi(z)|^{2}}{1-|z|^{2}},
\ee
(see\cite[p.198]{duren}) and from here
\be\label{eq6}
\left|\frac{1}{2}z^{3}\phi'(z)\right|\leq\frac{|z^{3}|}{2(1-|z|^{2})}(1-|\phi(z)|^{2})<1-|\phi(z)|^{2},
\ee
because $\frac{|z^{3}|}{2(1-|z|^{2})}<1$ for $|z|<r_{1}$.
Also,
\be\label{eq7}
|z^{2}\phi(z)|<|r_{1}^{2}\phi(z)|<\frac{1}{\sqrt{2}}|\phi(z)|,
\ee
since
$$r_{1}^{2}=0.7044\ldots<\frac{1}{\sqrt{2}}=0.7071\ldots.$$
Finally, by using \eqref{eq5},\eqref{eq6} and \eqref{eq7}, we have
\[
\begin{split}
\left|\arg \left[2\frac{1-\frac{1}{2}z^{3}\phi'(z)}{1+z^{2}\phi(z)}\right]\right|&\leq |\arg\left[1-\frac{1}{2}z^{3}\phi'(z)\right]|+|\arg(1+z^{2}\phi(z))|\\
&<\arcsin(1-|\phi(z)|^{2})+\arcsin\left(\frac{1}{\sqrt{2}}|\phi(z)|\right)\\
&=\arcsin \sqrt{1-\frac{1}{2}|\phi(z)|^{2}}\\
&\leq\arcsin1\\
&=\frac{\pi}{2},
\end{split}
\]
which implies ${\rm Re}[D(f;z)]>0$.

\medskip
 (ii) Since $f\in \mathcal{S}^{\star}(1/2)$, we can put
$$\frac{zf'(z)}{f(z)}=\frac{1}{1-\omega (z)},$$
where $\omega $ is analytic in $\ID$,  $\omega(0)=0$  and
$|\omega(z)|<1$ for all $z\in\ID$. From here we have that
$$\frac{zf''(z)}{f'(z)}=\frac{z\omega'(z)}{1-\omega (z)}+\frac{1}{1-\omega (z)}-1,$$
and so
$$D(f;z)=2-\frac{z\omega'(z)-\omega(z)}{1-\omega(z)}.$$
Since $\omega(0)=0$  and $|\omega(z)|<1$, $z\in\ID$, implies that $\left|\frac{\omega(z)}{z}\right|\leq 1$, $z\in\ID$, then by using the estimate  \eqref{eq5} (with $\frac{\omega(z)}{z}$ in stead of $\phi$), we obtain
\be\label{eq8}
|z\omega'(z)-\omega(z)|\leq\frac{r^{2}-|\omega(z)|^{2}}{1-r^{2}},
\ee
(where $|z|=r$ and
$|\omega(z)|\leq r$). Further, we have
\[
\begin{split}
{\rm Re}[D(f;z)]& \geq 2-\frac{|z\omega'(z)-\omega(z)|}{1-|\omega(z)|}\\
& \geq 2- \frac{1}{1-r^{2}}\frac{r^{2}-|\omega(z)|^{2}}{1-|\omega(z)|} \\
& = 2-\frac{1}{1-r^{2}}\varphi(t),
\end{split}
\]
where we put $|\omega(z)|=t,$ $0\leq t \leq r$ and $ \varphi(t)=\frac{r^{2}-t^{2}}{1-t}$. By elementary calculation we obtain that $\varphi(t)\leq 2(1-\sqrt{1-r^{2}})$ for $t\in [0,r].$ This  implies that
$${\rm Re}[D(f;z)]\geq 2-\frac{2(1-\sqrt{1-r^{2}})}{1-r^{2}}=2\frac{\sqrt{1-r^{2}}-r^{2}}{1-r^{2}}>0,$$
since $|z|=r<\sqrt{\frac{\sqrt{5}-1}{2}}=r_2$.

\medskip
(iii) For $f\in \mathcal{G}$ in \cite{JO_95} is proven that
$$\frac{zf'(z)}{f(z)}\prec \frac{1-z}{1-\frac{z}{2}},$$
i.e., that
$$\frac{zf'(z)}{f(z)}=\frac{1-\omega(z)}{1-\frac{\omega(z)}{2}},$$
where $\omega $ is analytic in $\ID$ such that $\omega(0)=0$  and
$|\omega(z)|<1$ for $z\in\ID$. From the last relation we easily obtain
$$D(f;z)=2-\frac{\omega(z)}{2-\omega(z)}+\frac{z}{(1-\omega(z))(2-\omega(z))}\omega '(z)$$
and from here
$${\rm Re}[D(f;z)]\geq 2-\frac{|\omega(z)|}{2-|\omega(z)|}-\frac{|z|}{(1-|\omega(z)|)(2-|\omega(z)|)}|\omega '(z)|.$$
Applying the inequality \eqref{eq5}, we give
$${\rm Re}[D(f;z)]\geq 2-\frac{|\omega(z)|}{2-|\omega(z)|}-\frac{|z|}{(1-|\omega(z)|)(2-|\omega(z)|)}\frac{1-|\omega(z)|^{2}}{1-|z|^{2}},$$
or, if we  use $|\omega(z)|\leq r$, where $|z|=r$:
$${\rm Re}[D(f;z)]\geq 2-\frac{r}{1-r}=\frac{2-3r}{1-r}>0,$$
since $r<\frac{2}{3}=r_3.$

\medskip
(iv) We can use relation
\eqref{eq5} and the same method as in the previous cases. Namely, now we can put
$$\frac{zf'(z)}{f(z)}=\frac{1+\omega(z)}{1-\omega (z)},$$
where $\omega $ is analytic in $\ID$, $\omega(0)=0$  and
$|\omega(z)|<1$ for $z\in\ID$. Then,
$$\frac{zf''(z)}{f'(z)}=2\frac{z\omega'(z)}{1-\omega ^{2}(z)}+\frac{1+\omega(z)}{1-\omega (z)}-1$$ and after that
$$D(f;z)=\frac{2}{1-\omega (z)}-2\frac{z\omega'(z)}{1-\omega^{2}(z)}.$$
Finally, we have (using $|\omega(z)|\leq r$, where $|z|=r$):
\[
\begin{split}
{\rm Re}[D(f;z)]& \ge {\rm Re}\frac{2}{1-\omega(z)} -2\frac{|z| |\omega'(z)|}{1-|\omega(z)|^{2}}\\
&\geq \frac{2}{1+|\omega(z)|}-2\frac{|z|}{1-|\omega(z)|^{2}}
\frac{1-|\omega(z)|^{2}}{1-|z|^{2}}\\
&\geq \frac{2}{1+r}-\frac{2r}{1-r^{2}}\\
&= 2\frac{1-2r}{1-r^{2}}>0,
\end{split}
\]
if $|z|=r<\frac{1}{2}=r_4.$

\medskip
(v) If $f\in \mathcal{S}$ then, from the classical result (see \cite[p.32]{duren}), we have
$$\left|\log \frac{zf'(z)}{f(z)}\right|\leq \log \frac{1+r}{1-r},\quad |z|=r<1  .$$
If we put $w=\log \frac{zf'(z)}{f(z)}$ and $R=\log \frac{1+r}{1-r}$, then we have
$\frac{zf'(z)}{f(z)}=e^{w}$, where $|w|\leq R $. If we choose $r\leq \text{tanh}\frac{1}{2}=\frac{e-1}{e+1}=0.46...$,
then we have $R\leq 1$. For such $R$ the function  $e^{w}$ is convex with positive real coefficients
that maps the unit disk onto a region that is symmetric with respect to the real axes ($\overline{e^w}=e^{\overline{w}}$), with diameter end points for $w=-1$ and $w=1$. This implies that
$${\rm Re}\frac{zf'(z)}{f(z)}={\rm Re}(e^{w})\geq e^{-R}=\frac{1-r}{1+r}.$$
Also, from the relation for functions from the class $\mathcal{S}$ (see \cite[Theorem 2.4, p.32]{duren}) we have
$$\left|\frac{zf''(z)}{f'(z)}-\frac{2r^{2}}{1-r^{2}}\right|\leq\frac{4r}{1-r^{2}}$$
and from here
$${\rm Re}\frac{zf''(z)}{f'(z)}\leq \frac{2r^{2}}{1-r^{2}}+\frac{4r}{1-r^{2}}=2\frac{2r+r^{2}}{1-r^{2}}.$$
Finally,
$${\rm Re}[D(f;z)]\geq 2\frac{1-r}{1+r}-2\frac{2r+r^{2}}{1-r^{2}}=2\frac{1-4r}{1-r^{2}}>0 $$
if $|z|=r<\frac{1}{4}=r_5$.
\end{proof}

\end{document}